\documentclass[12pt,a4paper]{article}

\usepackage{amsmath}
\usepackage{amssymb}
\usepackage{t1enc}
\usepackage[latin1]{inputenc}
\usepackage[english]{babel}
\usepackage{amsfonts}
\usepackage{latexsym}
\usepackage{bm}
\newtheorem{theorem}{Theorem}

\newtheorem{lemma}{Lemma}

\newtheorem{claim}{Claim}

\begin{document}

\title{The maximum sum of sizes of cross-intersecting families of subsets of a set}

\author{Peter Borg\\[5mm]
Department of Mathematics\\
Faculty of Science\\
University of Malta\\
Malta\\
\texttt{peter.borg@um.edu.mt} \\[5mm]
\and
Carl Feghali\\ [5mm]
Computer Science Institute of Charles University\\
Prague\\ 
Czech Republic\\
\texttt{feghali.carl@gmail.com}
}

\date{}
\maketitle

\begin{abstract}
A set of sets is called a \emph{family}. Two families $\mathcal{A}$ and $\mathcal{B}$ of sets are said to be \emph{cross-intersecting} if each member of $\mathcal{A}$ intersects each member of $\mathcal{B}$. For any two integers $n$ and $k$ with $1 \leq k \leq n$, let ${[n] \choose \leq k}$ denote the family of subsets of $[n] = \{1, \dots, n\}$ that have at most $k$ elements. We show that if $\mathcal{A}$ is a non-empty subfamily of ${[n] \choose \leq r}$, $\mathcal{B}$ is a non-empty subfamily of ${[n] \choose \leq s}$, $r \leq s$, and $\mathcal{A}$ and $\mathcal{B}$ are cross-intersecting, then 
\[|\mathcal{A}| + |\mathcal{B}| \leq 1 + \sum_{i=1}^s \left({n \choose i} - {n-r \choose i} \right),\] 
and equality holds if $\mathcal{A} = \{[r]\}$ and $\mathcal{B}$ is the family of sets in ${[n] \choose \leq s}$ that intersect $[r]$. 
\end{abstract}

\section{Introduction} \label{Intro}
Unless stated otherwise, we shall use small letters such as $x$ to
denote non-negative integers or elements of a set,
capital letters such as $X$ to denote sets, and calligraphic
letters such as $\mathcal{F}$ to denote \emph{families}
(that is, sets whose members are sets themselves). It is to be
assumed that arbitrary sets and families are finite. We
call a set $A$ an \emph{$r$-element set} if its size $|A|$ is $r$, that is, if it contains exactly $r$ elements.

The set $\{1, 2, \dots\}$ of positive integers is denoted by $\mathbb{N}$. For $m \geq 0$ and $n \geq 0$, the set $\{i \in \mathbb{N} \colon m \leq i \leq n\}$ is denoted by $[m, n]$. We abbreviate $[1,n]$ to $[n]$. Note that $[0]$ is the empty set $\emptyset$, and $[n] = \{1, \dots, n\}$ for $n \geq 1$. For a set $X$, the \emph{power set of $X$} (the set of subsets of $X$) is denoted by $2^X$. For any integer $r \geq 0$, the family of $r$-element subsets of $X$ is denoted by $X \choose r$, and the family of subsets of $X$ of size at most $r$ is denoted by $X \choose \leq r$. Thus, ${X \choose \leq r} = \bigcup_{i = 0}^r {X \choose i}$. If $x \in X$ and $\mathcal{F} \subseteq 2^X$, then we denote the family of sets in $\mathcal{F}$ which contain $x$ by $\mathcal{F}(x)$. We call $\mathcal{F}(x)$ a \emph{star of $\mathcal{F}$} if $\mathcal{F}(x) \neq \emptyset$.

We say that a set $A$ \emph{intersects} a set $B$ if $A$ and $B$ have at least one common element (that is, if $A \cap B \neq \emptyset$). A family $\mathcal{A}$ is said to be \emph{intersecting} if every two sets in $\mathcal{A}$ intersect. Note that the stars of a family $\mathcal{F}$ are the simplest intersecting subfamilies of $\mathcal{F}$. If $\mathcal{A}$ and $\mathcal{B}$ are families such that each set in $\mathcal{A}$ intersects each set in $\mathcal{B}$, then $\mathcal{A}$ and $\mathcal{B}$ are said to be \emph{cross-intersecting}. 

One of the most popular endeavours in extremal set theory is that of determining the size of a largest intersecting subfamily of a given family $\mathcal{F}$. This started in \cite{EKR}, which features the classical result, known as the Erd\H os-Ko-Rado (EKR) Theorem, that says that if $r \leq n/2$, then the size of a largest intersecting subfamily of ${[n] \choose r}$ is the size ${n-1 \choose r-1}$ of any star of ${[n] \choose r}$. There are many proofs of the EKR Theorem. Two of them are particularly short and beautiful: Katona's \cite{K}, which introduced the elegant cycle method, and Daykin's \cite{D}, which uses a fundamental result known as the Kruskal--Katona Theorem \cite{Ka,Kr} (see also \cite{FF,Kat}). The EKR Theorem gave rise to some of the highlights in extremal set theory \cite{AK1,F_t1,HM,Kat,W} and inspired many results, including generalizations (see, for example, \cite{Borg1,T}), that establish how large a system of sets can be under certain intersection conditions; see \cite{Borg7,DF,F2,F,FTsurvey,HST,HT}.

A natural question to ask about intersecting subfamilies of a given family $\mathcal{F}$ is how large they can be. For
cross-intersecting families, two natural parameters arise: the
sum and the product of sizes of the families. The problem of maximizing the sum or the product of sizes of cross-$t$-intersecting subfamilies of a given family $\mathcal{F}$ has been attracting much attention. 
Many of the results to date are referenced in \cite{Borg8,Borgmaxprod,BorgJLMS}.

Hilton and Milner \cite{HM} proved that, for $1 \leq r \leq n/2$, if
$\mathcal{A}$ and $\mathcal{B}$ are non-empty cross-intersecting subfamilies of ${[n] \choose r}$, then $|\mathcal{A}| +
|\mathcal{B}| \leq {n \choose r} - {n-r \choose r} + 1$, and equality holds if $\mathcal{A} = \{[r]\}$ and $\mathcal{B} = \{B \in {[n] \choose r} \colon B \cap [r] \neq \emptyset\}$. To the best of the authors' knowledge, this was the first result on the sizes of cross-intersecting families. 
Simpson \cite{Simpson} obtained a streamlined proof by means of the \emph{compression} (also known as \emph{shifting}) technique, which was introduced in the seminal EKR paper \cite{EKR} and has proven to be a very useful tool in extremal set theory (\cite{F} is a recommended survey on the properties and uses of compression operations). Frankl and Tokushige \cite{FT} instead used the Kruskal--Katona Theorem to establish the following stronger result: if $1 \leq r \leq s$, $n \geq r+s$, $\mathcal{A} \subseteq {[n] \choose r}$, $\mathcal{B} \subseteq {[n] \choose s}$, and $\mathcal{A}$ and $\mathcal{B}$ are cross-intersecting and non-empty, then $|\mathcal{A}| + |\mathcal{B}| \leq 1 + {n \choose s} - {n-r \choose s}$, and equality holds if $\mathcal{A} = \{[r]\}$ and $\mathcal{B} = \{B \in {[n] \choose s} \colon B \cap [r] \neq \emptyset\}$. The attainable upper bound on the maximum product of sizes for $1 \leq r \leq s \leq n/2$ was established in \cite{MT,Pyber}.

In this paper, we solve the analogous maximum sum problem for the case where $\mathcal{A} \subseteq {[n] \choose \leq r}$ and $\mathcal{B} \subseteq {[n] \choose \leq s}$, using the compression technique. The following is our result, proved in the next section. 

\begin{theorem}\label{mainthm} If $n \geq 1$, $1 \leq r \leq s$, $\mathcal{A} \subseteq {[n] \choose \leq r}$, $\mathcal{B} \subseteq {[n] \choose \leq s}$, and $\mathcal{A}$ and $\mathcal{B}$ are cross-intersecting and non-empty, then 
\[|\mathcal{A}| + |\mathcal{B}| \leq 1 + \sum_{i=1}^s \left({n \choose i} - {n-r \choose i} \right),\] 
and equality holds if $\mathcal{A} = \{[r]\}$ and $\mathcal{B} = \{B \in {[n] \choose \leq s} \colon B \cap [r] \neq \emptyset\}$.
\end{theorem}
The analogous problem for the product of sizes was solved in \cite{BorgBLMS}.

\section{Proof of Theorem~\ref{mainthm}} \label{Proof1}

We now prove Theorem~\ref{mainthm}. 

For any $i, j \in [n]$, let $\delta_{i,j}
\colon 2^{[n]} \rightarrow 2^{[n]}$ be defined by
\[ \delta_{i,j}(A) = \left\{ \begin{array}{ll}
(A \backslash \{j\}) \cup \{i\} & \mbox{if $j \in A$ and $i \notin
A$;}\\
A & \mbox{otherwise,}
\end{array} \right. \]
and let $\Delta_{i,j} \colon 2^{2^{[n]}} \rightarrow 2^{2^{[n]}}$ be the \emph{compression operation} defined by
\[\Delta_{i,j}(\mathcal{A}) = \{\delta_{i,j}(A) \colon A \in
\mathcal{A}\} \cup \{A \in
\mathcal{A} \colon \delta_{i,j}(A) \in \mathcal{A}\}.\]
Note that $\Delta_{i,j}$ preserves the size of a family $\mathcal{A}$, that is,
\[|\Delta_{i,j}(\mathcal{A})| = |\mathcal{A}|.\] 
We will need this equality together with the following basic fact, which we prove for completeness.

\begin{lemma}\label{compcross} If $\mathcal{A}$ and $\mathcal{B}$ are cross-intersecting subfamilies of $2^{[n]}$, then, for any $i, j \in [n]$, $\Delta_{i,j}(\mathcal{A})$ and $\Delta_{i,j}(\mathcal{B})$ are cross-intersecting subfamilies of $2^{[n]}$.
\end{lemma}
\textbf{Proof.} Suppose $A \in \Delta_{i,j}(\mathcal{A})$ and $B \in \Delta_{i,j}(\mathcal{B})$. If $A \in \mathcal{A}$ and $B \in \mathcal{B}$, then $A \cap B \neq \emptyset$. Suppose $A \notin \mathcal{A}$ or $B \notin \mathcal{B}$. We may assume that $A \notin \mathcal{A}$. Then, $A = \delta_{i,j}(A') \neq A'$ for some $A' \in \mathcal{A}$, so $i \notin A'$, $j \in A'$, $i \in A$, and $j \notin A$. Suppose $A \cap B = \emptyset$. Then, $i \notin B$, $B \in \mathcal{B} \cap \Delta_{i,j}(\mathcal{B})$, and hence $B, \delta_{i,j}(B) \in \mathcal{B}$. Thus, $A' \cap B \neq \emptyset$ and $A' \cap \delta_{i,j}(B) \neq \emptyset$. Since $A \cap B = \emptyset$ and $A' \cap B \neq \emptyset$, $A' \cap B = \{j\}$. This yields $A' \cap \delta_{i,j}(B) = \emptyset$, a contradiction.~\hfill{$\Box$}\\

If $i < j$, then $\Delta_{i,j}$ is called a \emph{left-compression}. A family $\mathcal{F} \subseteq 2^{[n]}$ is said to be \emph{compressed} if $\Delta_{i,j}(\mathcal{F}) = \mathcal{F}$ for every $i,j \in [n]$ with $i < j$ (that is, if $\mathcal{F}$ is invariant under left-compressions). Thus, $\mathcal{F}$ is compressed if and only if $(F \backslash \{j\}) \cup \{i\} \in \mathcal{F}$ for every $i, j \in [n]$ and every $F \in \mathcal{F}$ such that $i < j \in F$ and $i \notin F$. 

A subfamily $\mathcal{A}$ of $2^{[n]}$ that is not compressed can be transformed to a compressed subfamily of $2^{[n]}$ as follows. We choose one of the left-compressions that change $\mathcal{A}$, and we apply it to $\mathcal{A}$ to obtain a new subfamily of $2^{[n]}$. We keep on repeating this (always applying a left-compression to the last family obtained) until a family that is invariant under each left-compression is obtained (such a point is indeed reached, because if $\Delta_{i,j}(\mathcal{F}) \neq \mathcal{F} \subseteq 2^{[n]}$ and $i < j$, then $0 < \sum_{G \in \Delta_{i,j}(\mathcal{F})} \sum_{b \in G} b < \sum_{F \in \mathcal{F}} \sum_{a \in F} a$).

If $\mathcal{A}, \mathcal{B} \subseteq 2^{[n]}$ such that $\mathcal{A}$ and $\mathcal{B}$ are cross-intersecting, then, by Lemma~\ref{compcross}, we can obtain $\mathcal{A}^*, \mathcal{B}^* \subseteq 2^{[n]}$ such that $\mathcal{A}^*$ and $\mathcal{B}^*$ are compressed and cross-intersecting, $|\mathcal{A}^*| = |\mathcal{A}|$, and $|\mathcal{B}^*| = |\mathcal{B}|$. Indeed, similarly to the procedure above, if we can find a left-compression that changes at least one of $\mathcal{A}$ and $\mathcal{B}$, then we apply it to both $\mathcal{A}$ and $\mathcal{B}$, and we keep on repeating this (always performing this on the last two families obtained) until we obtain $\mathcal{A}^*, \mathcal{B}^* \subseteq 2^{[n]}$ such that $\mathcal{A}^*$ and $\mathcal{B}^*$ are compressed.
\\
\\
\textbf{Proof of Theorem~\ref{mainthm}.} We prove the result by induction on $n$. The result is trivial for $n \leq 2$. Consider $n \geq 3$.  

Suppose $r \geq n$. Then, $s \geq n$ and ${[n] \choose \leq r} = {[n] \choose \leq s} = 2^{[n]}$. Since $\mathcal{A}$ and $\mathcal{B}$ are cross-intersecting, $[n] \backslash A \notin \mathcal{B}$ for each $A \in \mathcal{A}$. Thus, $\mathcal{B} \subseteq 2^{[n]} \backslash \{[n] \backslash A \colon A \in \mathcal{A}\}$, and hence $|\mathcal{B}| \leq 2^n - |\mathcal{A}|$. We have $|\mathcal{A}| + |\mathcal{B}| \leq 2^n = 1 + |\{B \in 2^{[n]} \colon B \cap [n] \neq \emptyset\}| = 1 + |\{B \in {[n] \choose \leq s} \colon B \cap [r] \neq \emptyset\}|$. 

Now suppose $r < n$. If $s > n$, then ${[n] \choose \leq s} = 2^{[n]} = {[n] \choose \leq n}$. Thus, we may assume that $s \leq n$. Since $r < n$ and $r \leq s$, we have
\begin{equation} \mbox{$r < s$ \quad or \quad $s < n$.}  \label{nontrivialcase}
\end{equation}

As explained above, we may assume that $\mathcal{A}$ and $\mathcal{B}$ are compressed.

Let $\mathcal{A}_0 = \{A \in \mathcal{A} \colon n \notin A\}$, $\mathcal{A}_1 = \{A \backslash \{n\} \colon n \in A \in
\mathcal{A}\}$, $\mathcal{B}_0 = \{B \in \mathcal{B} \colon n \notin B\}$, and $\mathcal{B}_1 = \{B \backslash \{n\} \colon n \in B \in
\mathcal{B}\}$. We have $\mathcal{A}_0 \subseteq {[n-1] \choose \leq r}$, $\mathcal{A}_1 \subseteq {[n-1] \choose \leq r-1}$, $\mathcal{B}_0 \subseteq {[n-1] \choose \leq s}$, and $\mathcal{B}_1 \subseteq {[n-1] \choose \leq s-1}$. Clearly, $\mathcal{A}_0$ and $\mathcal{B}_0$ are cross-intersecting. Since $\mathcal{A}$ and $\mathcal{B}$ are compressed, we clearly have that $\mathcal{A}_0$, $\mathcal{A}_1$, $\mathcal{B}_0$, and $\mathcal{B}_1$ are compressed. Thus, $[r'] \in \mathcal{A}_0$ for some $r' \in [r]$, and if $s < n$, then $[s'] \in \mathcal{B}_0$ for some $s' \in [s]$.

Let $\mathcal{C} = \{A \in \mathcal{A}_1 \colon A \cap B = \emptyset \mbox{ for some } B \in \mathcal{B}_1\}$. For each $C \in \mathcal{C}$, let $\bar{C} = [n-1] \backslash C$, $C' = C \cup \{n\}$, and $\bar{C}' = \bar{C} \cup \{n\}$. Let $\bar{\mathcal{C}} = \{\bar{C} \colon C \in \mathcal{C}\}$. For each $C \in \mathcal{C}$, $C' \in \mathcal{A}$ as $C \in \mathcal{A}_1$.

Suppose $C \in \mathcal{C}$. Let $\mathcal{D}_C = \{B \in \mathcal{B}_1 \colon B \cap C = \emptyset\}$. Suppose that there exists some $B \in \mathcal{D}_C$ such that $B \neq \bar{C}$. Then, $B \subsetneq [n-1] \backslash C$, and hence $[n-1] \backslash (B \cup C) \neq \emptyset$. Let $x \in [n-1] \backslash (B \cup C)$. Since $B \in \mathcal{B}_1$, $B \cup \{n\} \in \mathcal{B}$. Let $D = \delta_{x,n}(B \cup \{n\})$. Since $x \notin B \cup \{n\}$, $D = B \cup \{x\}$. Since $\mathcal{B}$ is compressed, $D \in \mathcal{B}$. However, since $x \notin C'$ and $B \cap C = \emptyset$, we have $C' \cap D = \emptyset$, which is a contradiction as $\mathcal{A}$ and $\mathcal{B}$ are cross-intersecting. Thus, $\mathcal{D}_C \subseteq \{\bar{C}\}$. Since $C \in \mathcal{C}$, we have $\mathcal{D}_C \neq \emptyset$, so $\mathcal{D}_C = \{\bar{C}\}$.

We have therefore shown that for any $C \in \mathcal{C}$, $\bar{C}$ is the unique set in $\mathcal{B}_1$ that does not intersect $C$. 
%
%
Since $C' \in \mathcal{A}$ and $C' \cap \bar{C} = \emptyset$, the cross-intersection condition gives us $\bar{C} \notin \mathcal{B}_0$. Thus,
\begin{equation} \mathcal{B}_0 \cap \bar{\mathcal{C}}  = \emptyset. \label{2.5}
\end{equation}

Let $\mathcal{A}_1' = \mathcal{A}_1 \backslash \mathcal{C}$ and $\mathcal{B}_0' = \mathcal{B}_0 \cup \bar{\mathcal{C}}$. Clearly, $\mathcal{A}_1' \subseteq {[n-1] \choose \leq r-1}$, $\mathcal{B}_0' \subseteq {[n-1] \choose \leq s}$, $\mathcal{A}_0$ and $\mathcal{B}_0'$ are cross-intersecting (as $\bar{\mathcal{C}} \subseteq \mathcal{B}_1$), and $\mathcal{A}_1'$ and $\mathcal{B}_1$ are cross-intersecting. 

\begin{claim} \label{claim1} $|\mathcal{A}_0| + |\mathcal{B}_0'| \leq 1 + |\{B \in {[n-1] \choose \leq s} \colon B \cap [r] \neq \emptyset\}|$. 
\end{claim}
\textbf{Proof.} Since $\mathcal{A}_0$ is non-empty (as $[r'] \in \mathcal{A}_0$), the claim follows by the induction hypothesis if $\mathcal{B}_0'$ is non-empty too, and this is the case if $s < n$ (as we then have $[s'] \in \mathcal{B}_0'$). Suppose $s = n$ and $\mathcal{B}_0' = \emptyset$. By (\ref{nontrivialcase}), $r \leq s-1$. Since $s = n$, the sets in $\mathcal{A}_0$ intersect $[s-1]$ (note that $\emptyset \notin \mathcal{A}$ as $\mathcal{B} \neq \emptyset$), so $\mathcal{A}_0$ and $\{[s-1]\}$ are cross-intersecting. By the induction hypothesis, $|\mathcal{A}_0| + |\{[s-1]\}| \leq 1 + |\{B \in {[n-1] \choose \leq s-1} \colon B \cap [r] \neq \emptyset\}|$, so $|\mathcal{A}_0| + |\mathcal{B}_0'| = |\mathcal{A}_0| \leq |\{B \in {[n-1] \choose \leq s} \colon B \cap [r] \neq \emptyset\}|$. \hfill{$\Box$}

\begin{claim} \label{claim2} $|\mathcal{A}_1'| + |\mathcal{B}_1| \leq |\{B \in {[n-1] \choose \leq s-1} \colon B \cap [r] \neq \emptyset\}|$.
\end{claim}
\textbf{Proof.} The claim is immediate if $\mathcal{A}_1' = \emptyset$ and $\mathcal{B}_1 = \emptyset$. 

If $\mathcal{A}_1' \neq \emptyset$ and $\mathcal{B}_1 \neq \emptyset$, then, by the induction hypothesis, $|\mathcal{A}_1'| + |\mathcal{B}_1| \leq 1 + |\{B \in {[n-1] \choose \leq s-1} \colon B \cap [r-1] \neq \emptyset\}| \leq |\{B \in {[n-1] \choose \leq s-1} \colon B \cap [r] \neq \emptyset\}|$.

If $\mathcal{A}_1' = \emptyset$ and $\mathcal{B}_1 \neq \emptyset$, then $|\mathcal{A}_1'| + |\mathcal{B}_1| = |\mathcal{B}_1| \leq |\{B \in {[n-1] \choose \leq s-1} \colon B \cap [r'] \neq \emptyset\}|$ (as $[r'] \in \mathcal{A}$), so $|\mathcal{A}_1'| + |\mathcal{B}_1| \leq |\{B \in {[n-1] \choose \leq s-1} \colon B \cap [r] \neq \emptyset\}|$.

Finally, suppose that $\mathcal{A}_1' \neq \emptyset$ and $\mathcal{B}_1 = \emptyset$.

Suppose $\mathcal{B} = \{[n]\}$. Then, $s = n$. By the definition of $\mathcal{A}_1'$, we have $\emptyset \notin \mathcal{A}_1'$, so $\mathcal{A}_1'$ and $\{[s-1]\}$ are cross-intersecting (as $s = n$). By the induction hypothesis, $|\mathcal{A}_1'| + |\{[s-1]\}| \leq 1 + |\{B \in {[n-1] \choose \leq s-1} \colon B \cap [r-1] \neq \emptyset\}|$, so $|\mathcal{A}_1'| + |\mathcal{B}_1| = |\mathcal{A}_1'| \leq |\{B \in {[n-1] \choose \leq s-1} \colon B \cap [r] \neq \emptyset\}|$.

Now suppose $\mathcal{B} \neq \{[n]\}$. Then, since $\mathcal{B}$ is compressed and non-empty, $[s^*] \in \mathcal{B}$ for some $s^* \in [s] \cap [n-1]$. Thus, $\mathcal{A}_1'$ and $\{[s^*]\}$ are cross-intersecting. If $s^* \leq s-1$, then the claim follows as in the preceding paragraph. Suppose $s^* = s$. Let $\mathcal{E} = \{A \in \mathcal{A}_1' \colon 1 \in A\}$ and $\mathcal{E}' = \mathcal{A}_1' \backslash \mathcal{E}$. Then, $\mathcal{E}'$ is a subfamily of ${[2,n-1] \choose \leq r-1}$ and its sets intersect the $(s-1)$-element set $[2,s]$. Let $\mathcal{F} = \{B \in {[n-1] \choose \leq s-1} \colon 1 \in B\}$ and $\mathcal{F}' = \{B \in {[2,n-1] \choose \leq s-1} \colon B \cap [2,r] \neq \emptyset\}$. Since $\mathcal{A}_1' \subseteq {[n-1] \choose \leq r-1} \subseteq {[n-1] \choose \leq s-1}$, $|\mathcal{E}| \leq |\mathcal{F}|$. By the induction hypothesis, $|\mathcal{E}'| + |\{[2,s]\}| \leq |\{[2,r]\}| + |\mathcal{F}'|$. We have $|\mathcal{A}_1'| + |\mathcal{B}_1| = |\mathcal{A}_1'| = |\mathcal{E}| + |\mathcal{E}'| \leq |\mathcal{F}| + |\mathcal{F}'| = |\{B \in {[n-1] \choose \leq s-1} \colon B \cap [r] \neq \emptyset\}|$. \hfill{$\Box$}
\\

We have
\begin{align} |\mathcal{A}| + |\mathcal{B}| &= |\mathcal{A}_0| + |\mathcal{A}_1| + |\mathcal{B}_0| + |\mathcal{B}_1| \nonumber \\
&= (|\mathcal{A}_0| + |\mathcal{B}_0'|) + (|\mathcal{A}_1'| + |\mathcal{B}_1|) + |\mathcal{C}| - |\bar{\mathcal{C}}| \quad \mbox{(by (\ref{2.5}))} \nonumber \\
&= (|\mathcal{A}_0| + |\mathcal{B}_0'|) + (|\mathcal{A}_1'| + |\mathcal{B}_1|). \nonumber
\end{align}
Therefore, by Claims~\ref{claim1} and \ref{claim2}, $|\mathcal{A}| + |\mathcal{B}| \leq 1 + |\{B \in {[n-1] \choose \leq s} \colon B \cap [r] \neq \emptyset\}| + |\{B \in {[n-1] \choose \leq s-1} \colon B \cap [r] \neq \emptyset\}| = 1 + |\{B \in {[n] \choose \leq s} \colon B \cap [r] \neq \emptyset\}|$. \hfill{$\Box$}
\\
\\
\textbf{Acknowledgements.} Peter Borg was supported by grant MATRP14-20 of the University of Malta. Carl Feghali was supported by grant 19-21082S of the Czech Science Foundation. 

\end{document}